\documentclass[12pt]{amsart}
\usepackage{amsmath}
\usepackage{graphicx}
\usepackage{latexsym}
\usepackage{pinlabel}
\usepackage{verbatim}
\usepackage[all]{xy}
\usepackage{url}

\newtheorem{thm}{Theorem}
\newtheorem*{mthm}{Main Theorem}
\newtheorem*{fscriterion}{Fiber Sum Criterion}
\newtheorem*{sscriterion}{Section Sum Criterion}

\newtheorem{theorem}[thm]{Theorem}

\theoremstyle{definition}

\newcommand{\R}{\mathbb{R}}
\newcommand{\Z}{\mathbb{Z}}

\def \x {\times}

\newcommand{\nc}{\newcommand}
\nc{\dmo}{\DeclareMathOperator}

\dmo{\MCG}{Mod}
\dmo{\Diff}{Diff}

\nc{\margin}[1]{\marginpar{\scriptsize #1}}

\begin{document}

\title[Indecomposable surface bundles]
{Indecomposable surface bundles \\ over surfaces}

\author[R. \.{I}. Baykur]{R. \.{I}nan\c{c} Baykur}
\address{Max Planck Institute for Mathematics, Bonn, Germany \newline
\indent Department of Mathematics, Brandeis University, Waltham, MA, USA}
\email{baykur@mpim-bonn.mpg.de, baykur@brandeis.edu}

\author[Dan Margalit]{Dan Margalit}
\address{School of Mathematics, Georgia Institute of Technology, Atlanta, GA, USA}
\email{margalit@math.gatech.edu}

\thanks{The first author was partially supported by the NSF grant DMS-0906912.  The second author was partially supported by an NSF CAREER grant and a fellowship from the Sloan Foundation.}

\begin{abstract}
For each pair of integers $g \geq 2$ and $h \geq 1$, we explicitly construct infinitely many fiber sum and section sum indecomposable genus $g$ surface bundles over genus $h$ surfaces whose total spaces are pairwise homotopy inequivalent.
\end{abstract}

\maketitle

\setcounter{secnumdepth}{2}
\setcounter{section}{0}

\section{Introduction}

A surface bundle over a surface is a surjective submersion $f\colon X \to B$, where $X$ and $B$ are closed, smooth, oriented $4$- and $2$-dimensional manifolds, respectively.  We say that $(X,f)$ is a genus $g$ surface bundle over a genus $h$ surface if the genus of a fiber $F$ is $g$ and the genus of the base $B$ is $h$.

There are two common ways to construct new surface bundles over surfaces from old ones, namely, by summing along fibers or along sections.  We now explain both constructions. 

First, suppose $(X_1, f_1)$ and $(X_2,f_2)$ are genus $g$ surface bundles over surfaces of genus $h_1$ and $h_2$, respectively.  Let $F_i$ be a fiber of $f_i$ for $i=1,2$. The \textit{fiber sum} of $(X_1, f_1)$ and $(X_2,f_2)$ is obtained by removing a fibered tubular neighborhood of each $F_i$ and then identifying the resulting boundaries via any fiber-preserving, orientation-reversing diffeomorphism. The end result is a genus $g$ bundle over a genus $h=h_1+h_2$ surface.

On the other hand, if $(X_1, f_1)$ and $(X_2,f_2)$ are genus $g_i$ surface bundles over surfaces of genus $h$, and if each $(X_i, f_i)$ has a section $S_i$ with self-intersection number $n_i$ so that $n_1 = - n_2$, we can take the \textit{section sum} of $(X_1, f_1)$ and $(X_2,f_2)$ along the $S_i$ in order to obtain a genus $g=g_1+g_2$ surface bundle over a genus $h$ surface, this time by removing a fibered tubular neighborhood of each $S_i$ and then identifying the resulting boundaries via any base-preserving, orientation-reversing diffeomorphism.

We say that a fiber sum is trivial if one of the two bundles is \linebreak $\Sigma \x S^2$.  Similarly, a section sum is trivial if one of the bundles is an $S^2$-bundle over a surface $\Sigma$.  If a surface bundle over a surface cannot be expressed as a nontrivial fiber/section sum, then we say it is \textit{fiber/section sum indecomposable}. 

\begin{mthm} \label{mainthm}
For any fixed $g \geq 2$ and $h \geq 1$, there are infinitely many fiber sum and section sum indecomposable genus $g$ surface bundles over genus $h$ surfaces whose total spaces are pairwise homotopy inequivalent.
\end{mthm}

The Main Theorem is proven by explicitly constructing surface bundles with the stated properties.  Consider the case $h=2$.  For any integer $n$, let $X_n \to \Sigma_2$ be the surface bundle prescribed by the monodromy factorization: 
\[
 [
 T_{c_2}^{-n}T_{c_1}^{n}
  \, , 
  W T_{c_5}^{n} W T_{c_5}^{-n}
 ]
\,
[
T_{c_3}^nT_{c_5}^{n}WT_{c_5}^{-1}
\, , 
T_{c_5}^nT_{c_1}^{-n}T_{c_5}^{-n}T_{c_1}^{n}
] 
 =1.
\]
where $W = T_{c_4}^{n} T_{c_5}^{-n} \cdots T_{c_{2g-1}}^{-n} T_{c_{2g}}^{n} T_{c_{2g-1}}^{-n} \cdots T_{c_5}^{-n}T_{c_4}^{n}$, and where the simple closed curves $c_1,\dots,c_{2g+1}$ are as shown in Figures~\ref{figure:braid} and~\ref{figure:braid2} and $T_{c_i}$ is the (right) Dehn twist about $c_i$.  By varying $n$ over the set of odd primes, we obtain the surface bundles promised by the Main Theorem.

Our bundles for $h \geq 3$ are pullbacks of the above bundles via a fixed covering map $\Sigma_h \to \Sigma_2$; see Step 7 in Section~\ref{sec:construction} for the explicit monodromy factorizations.

The hypotheses on genus in our Main Theorem are in fact necessary: When $g < 2$, the identity component of the group of orientation-preserving self-diffeomorphisms of the fiber is not simply-connected. So using twisted gluings one can decompose any genus zero or genus one bundle as a nontrivial fiber sum where one of the summands is a nontrivial ruled or elliptic fibration over the $2$-sphere. On the other hand, if $h=0$, one can see that there is a unique surface bundle of genus $g \geq 2$, namely the trivial bundle $\Sigma_g \x S^2$, which can be expressed as the section sum of trivial surface bundles of smaller fiber genera. 

In some cases, one can argue the existence of fiber sum or section sum indecomposable surface bundles using the topology of the underlying $4$-manifolds. For example, the surface bundles over genus two surfaces with nonzero signatures produced by Bryan and Donagi \cite{BD} are clearly fiber sum indecomposable because any surface bundle over a genus zero or one surface has signature zero \cite{AS, EKKOS}, and thus by Novikov additivity their fiber sum has signature zero \cite{Novikov}. Similarly, the genus three surface bundles of Endo, Korkmaz, Kotschick, Ozbagci, and Stipsicz  \cite{EKKOS} with nonzero signatures can be seen to be section sum indecomposable, since any surface bundle of fiber genus less than three is hyperelliptic, and hence has vanishing signature. However, it is a priori unclear that for either one of these families of bundles we can strike both indecomposability properties in question, nor can this signature argument can be employed to cover higher fiber or base genus examples.

{\bf Outline.} In Section~\ref{Sec: Algebraic} we will provide algebraic characterizations of fiber sum indecomposability and section sum indecomposability, respectively. The proof of Theorem~\ref{mainthm} will then rely on the advances in geometric group theory on embeddings of surface groups into mapping class groups, which we will review in Section~\ref{Sec: Embeddings}. We will focus on the simplest and most explicit embeddings of this sort, which factor through right angled Artin groups and braid groups. Section~\ref{Sec: Construction} is where we will prove Theorem~\ref{mainthm}, and discuss the underlying geometric structures. Unlike the examples discussed above, the total spaces of the bundles we construct will all have signature zero.

\bigskip {\bf Acknowledgments.} This paper was partly inspired by a question of Ursula Hamenst\"adt. We would like to thank her as well as John Etnyre, Sang-hyun Kim, and Michael L\"onne for helpful conversations.  We are especially grateful to Sang-hyun Kim for pointing out a mistake in an earlier draft.

\section{Algebraic characterizations of indecomposability} \label{Sec: Algebraic}

Let $\Sigma$ denote a compact, connected oriented surface with a finite set of marked points in its interior.  The \emph{mapping class group} of $\Sigma$, denoted $\MCG(\Sigma)$, is the group of isotopy classes of orientation-preserving self-diffeomorphisms of $\Sigma$ that fix $\partial\Sigma$ pointwise and preserve the set of marked points.

We denote by $\Sigma_g$ the closed, connected, orientable surface of genus $g$.  We also denote by $\Sigma_{g,1}$ and $\Sigma_g^1$ the surfaces obtained from $\Sigma_g$ by marking one point and by deleting the interior of an  embedded disk, respectively.  There are surjective homomorphisms $\MCG(\Sigma_g^1) \to \MCG(\Sigma_{g,1})$ and $\MCG(\Sigma_{g,1}) \to \MCG(\Sigma_g)$ obtained by collapsing the boundary to a marked point and by forgetting the marked point, respectively.

Say that $g \geq 2$, and let $B$ be any Hausdorff, paracompact space.  A classical result of Earle and Eells states that the connected components of the diffeomorphism group $\Diff(\Sigma_g)$ are contractible \cite{EE}.  It follows that there is a bijective correspondence:
\begin{eqnarray*}
\left\{\begin{array}{c}
\mbox{Genus $g$ surface}\\
\mbox{bundles over $B$}\\
\mbox{up to isomorphism}\\
\end{array}\right\}
&
\longleftrightarrow
&
\left\{\begin{array}{c}
\mbox{Homomorphisms}\\
\pi_1(B) \to \MCG(\Sigma_g) \\
\mbox{up to conjugacy}\\
\end{array}\right\}
\end{eqnarray*}
The homomorphism $\mu: \pi_1(B) \to \MCG(\Sigma_g)$ corresponding to a given bundle is called the \emph{monodromy} of the bundle.

Let $f : X \to \Sigma_h$ be a genus $g$ surface bundle over a genus $h$ surface. 

Choose generators $\alpha_j, \beta_j$ of $\pi_1(\Sigma_h)$ so that
\[ \prod_{j=1}^h [\alpha_j, \beta_j]=1.\]
Since this is the only defining relation for $\pi_1(\Sigma_h)$, a bundle $X \to \Sigma_h$ is completely determined by the images of the $\alpha_j$ and $\beta_j$ under the monodromy $\mu$.  In other words, genus $g$ bundles $f : X \to \Sigma_h$ are completely determined by choices of $\mu(\alpha_j),\mu(\beta_j) \in \MCG(\Sigma_g)$ satisfying the relation
\[ \prod_{j=1}^h [\mu(\alpha_j), \mu(\beta_j)]=1.\]
Such an expression is called a \emph{monodromy factorization} for $(X,f)$. 

The purpose of this section is to give algebraic interpretations of the fiber sum and section sum operations in terms of the corresponding monodromies.

\subsection{Fiber sums via free products}\label{section:fscrit}  Let $(X_1,f_1)$ and $(X_2,f_2)$ be genus $g$ surface bundles over surfaces of genus $h_1$ and $h_2$, respectively.  Let $\mu_i \colon \pi_1(B_i) \to \MCG(\Sigma_g)$ be the monodromy of the bundle $f_i$, for $i=1,2$.  There is an induced homomorphism
\[ \mu_1 \ast \mu_2 : \pi_1 (B_1) \ast \pi_1 (B_2) \to \MCG(\Sigma_g). \]
Let $B$ be a surface of genus $h=h_1+h_2$ obtained by taking the connected sum of the $B_i$, and let $\gamma$ denote the simple closed curve in $B$ along which the $B_i \setminus D^2$ are glued.  Base $\pi_1(B)$ at a point of $\gamma$.  There is a homomorphism
\[ \pi_1(B) \to \pi_1(B_1) \ast \pi_1(B_2) \]
induced by collapsing $\gamma$ to a point. The monodromy of the fiber sum of $(X_1,f_1)$ and $(X_2,f_2)$ is precisely the one induced by postcomposing the above map with $\mu_1 \ast \mu_2$.

We conclude that a surface bundle over a surface is fiber sum indecomposable if and only if its monodromy does not decompose into a nontrivial free product of surface bundle monodromies as above.  More precisely:

\begin{fscriterion}
A surface bundle over a surface is fiber sum decomposable if and only if the kernel of the monodromy contains a nontrivial separating simple closed curve.
\end{fscriterion}

We will make use of the following immediate consequence of the Fiber Sum Criterion:
\begin{quote}
\emph{If a surface bundle over a surface has injective monodromy, then the bundle is fiber sum indecomposable.}
\end{quote}

We can restate the Fiber Sum Criterion in terms of the monodromy factorization: a genus $g$ bundle over $B \cong \Sigma_h$ with monodromy $\mu$ is fiber sum decomposable if and only if there is a choice of generators $\alpha_j,\beta_j$ for $\pi_1(B)$, and a $1 \leq k < h$ so that
\[ [\mu(\alpha_1),\mu(\beta_1)]\cdots[\mu(\alpha_k),\mu(\beta_k)] = 1 \]
in $\MCG(\Sigma_g)$.

The fiber sum operation involves a choice of gluing map, and we can say precisely how this affects the resulting monodromy factorization: the different choices of gluing map correspond to all possible monodromy factorizations of the form
\begin{align*} \phi([\mu(\alpha_1)&,\mu(\beta_1)]\cdots[\mu(\alpha_k),\mu(\beta_k)])\phi^{-1} \\  &[\mu(\alpha_{k+1}),\mu(\beta_{k+1})] \cdots [\mu(\alpha_g),\mu(\beta_g)] = 1 \end{align*}
for $\phi \in \MCG(\Sigma_g)$.

\subsection{Section sums via direct products}\label{section:sscrit} Let $(X_i,f_i)$ be genus $g_i$ surface bundles over $\Sigma_h$, and let $S_i$ be a section of $f_i$ with self-intersection number $n_i$. Assume that $n_1 = - n_2$. We again denote the monodromy of $f_i$ by $\mu_i$, for $i=1,2$. 

The boundary of a fibered regular neighborhood of $S_i$ is a circle bundle over $\Sigma_h$ whose euler class is $n_i$.   Thus, there is a base-preserving, orientation reversing diffeomorphism between these two circle bundles. 

Let $F$ be a genus $g=g_1+g_2$ fiber obtained by taking the section sum of $(X_1,f_1)$ and $(X_2,f_2)$ along the $S_i$.  If $F_1$ and $F_2$ are fibers of $f_i$ over the same point, then $F$ is the connect sum of $F_1$ and $F_2$ along a simple closed curve $\gamma \subset F$.

By construction, the monodromy $\mu$ of the section sum preserves the simple closed curve $\gamma$ with orientation.  There is a short exact sequence
\[ 1 \to \langle T_\gamma \rangle \to \MCG(F,\vec \gamma) \to \MCG(\Sigma_{g_1,1}) \times \MCG(\Sigma_{g_1,1}) \to 1 \]
where $\MCG(F,\vec \gamma)$ is the subgroup of $\MCG(F)$ consisting of all elements that preserve the isotopy class of $\gamma$ with orientation \cite[Proposition 3.20]{FM}.  The monodromy $\mu$ is a lift of $\mu_1 \times \mu_2$ to $\MCG(F,\vec \gamma)$.

Conversely, if the monodromy of a surface bundle preserves the orientation of an oriented, essential, separating simple closed curve, then we can reverse the above process, and so the bundle is section sum decomposable.

\begin{sscriterion}
A surface bundle over a surface with fiber $F \cong \Sigma_g$ is section sum decomposable if and only if there is an oriented essential separating simple closed curve $\gamma \subset F$ whose isotopy class is preserved by each element of the monodromy.
\end{sscriterion}

An element of $\MCG(\Sigma_g)$ is \emph{irreducible} if it does not preserve the isotopy class of any homotopically essential 1-submanifold of $\Sigma_g$.  Similarly, a subgroup of $\MCG(\Sigma_g)$ is irreducible if there is no homotopically essential 1-submanifold preserved up to isotopy by each element of the group, and a monodromy $\pi_1(\Sigma_h) \to \MCG(\Sigma_g)$ is irreducible if its image is an irreducible subgroup.

We will use the following immediate consequence of the Section Sum Criterion:
\begin{quote}
\emph{If the monodromy of a surface bundle over a surface is irreducible, then the bundle is section sum indecomposable.}
\end{quote}
To show that a monodromy is irreducible, it is enough to find one irreducible element in its image.  In fact, the converse is also true: every infinite irreducible subgroup contains an irreducible element \cite[Corollary 7.14]{ivanov}.

We can again restate our criterion in terms of monodromy factorizations.  Let $f:X \to B$ be a genus $g$ surface bundle over a genus $h$ surface.  Fix a generating set $\alpha_i,\beta_i$ for $B$ with $[\alpha_1,\beta_1]\cdots[\alpha_h,\beta_h]=1$.  Fix a fiber $F$ and let $\mu : \pi_1(B) \to \MCG(F)$ be the monodromy. The bundle $(X,f)$ is section sum decomposable if and only if there is an oriented, essential, separating simple closed curve $\gamma \subset F$ with the property that, for each $\delta \in \{\alpha_i,\beta_i\}$, we have 
\[ \mu(\delta)(\gamma)=\gamma. \]
As such, we can write
\[ \mu(\delta) = \mu^+(\delta)\mu^-(\delta), \]
where $\mu^+(\delta)$ and $\mu^-(\delta)$ are supported on the subsurfaces to the left and right of $\gamma$, respectively.  It follows that
\[ \prod_{j=1}^h [\mu^+(\alpha_j),\mu^+(\beta_j)] = T_{\gamma}^n \ \ \text{and} \ \
\prod_{j=1}^h [\mu^-(\alpha_j),\mu^-(\beta_j)] = T_{\gamma}^{-n}. \]

So in this case the monodromy factorization can be written as 
\[ \prod_{j=1}^h [\mu^+(\alpha_j) \mu^-(\alpha_j) ,\mu^+(\beta_j) \mu^-(\beta_j)] = T_{\gamma}^{n} T_{\gamma}^{-n} = 1 .
\]

As in the fiber sum case, different choices of gluing maps give rise to different monodromies.  Given one monodromy $\mu$ as above, all other monodromies take the form:
\[ \mu(\delta) = \phi\mu^+(\delta)\phi^{-1}\mu^-(\delta), \]
for varying $\phi$ in $\MCG(\Sigma_{g_1}^1)$, the mapping class group of the closed subsurface to the left of $\gamma$ (as above, $\delta \in \{\alpha_i,\beta_i\}$); here $\phi$ is independent of $\delta$.

\section{From surface groups to mapping class groups} \label{Sec: Embeddings}

By the Fiber Sum Criterion and the Section Sum Criterion of Section~\ref{Sec: Algebraic}, we can construct indecomposable surface bundles over surfaces if we can construct injective irreducible homomorphisms from surface groups to mapping class groups.


The homomorphisms $\pi_1(\Sigma_h) \to \MCG(\Sigma_g)$ we construct in Section~\ref{Sec: Construction} will pass through the theories of right-angled Artin groups and braid groups:
\begin{quote}
\emph{Surface groups $\to$ Right-angled Artin groups $\to$ \newline \hspace*{1.15in} Braid groups $\to$ Mapping class groups}
\end{quote}
The goal of this section is to explain each of these relationships.  Before we begin in earnest, we recall the definition of a right-angled Artin group.

A \emph{right-angled Artin group} is a group defined by a presentation with generators $\{v_1,\dots,v_n\}$, and where each defining relation has the form $v_iv_j=v_jv_i$.  To any finite graph $\Gamma$, we can associate a right-angled Artin group $A(\Gamma)$ by taking one generator $v_i$ for each vertex and a relation $v_iv_j=v_jv_i$ whenever the corresponding vertices are connected by an edge.

\subsection{Surface groups into right-angled Artin groups}\label{sec:crisp wiest} We now give a method due to Crisp--Wiest for explicitly embedding surface groups into right-angled Artin groups \cite{CW}.

Let $\Gamma$ be a finite graph with vertex set $v_i$.  A \emph{$\Gamma$-dissection} of a closed, oriented surface $\Sigma$ is a finite union of labeled, oriented simple closed curves $\{c_i\}$ in $\Sigma$ with the following properties:
\begin{enumerate}
\item The label of each $c_i$ is an element of the set $\{v_i\}$ (multiple curves can have the same label).
\item If two curves $c_i$ and $c_j$ intersect then the corresponding vertices of $\Gamma$ are distinct and are connected by an edge.
\item The complement of $\cup c_i$ in $\Sigma$ is a disjoint union of disks.
\end{enumerate}

Let $\gamma$ be any oriented loop in $\Sigma_h$ that does not pass through any intersection points of the curves $v_i$.  We can read off a word $w(\gamma)$ in the $v_i$ and their inverses by keeping track of the order and the directions (right-to-left or left-to-right) in which $\gamma$ crosses the $v_i$-curves.

The rule $\gamma \mapsto w(\gamma)$ gives rise to a homomorphism $\pi_1(\Sigma) \to A(\Gamma)$ called the \emph{label-reading map}.   To see that the label-reading map is well-defined, notice that, if we pass $\gamma$ across an intersection point of curves labeled $v_i$ and $v_j$, then that corresponds to swapping $v_i$ and $v_j$ in $w(\gamma)$.  By part 2 of the definition of a $\Gamma$-dissection and the defining relations of $A(\Gamma)$, the resulting element of $A(\Gamma)$ is the same as before.  Any two representatives of an element of $\pi_1(\Sigma_h)$ differ by such moves, and so well-definedness follows.

We will be interested in when the label-reading map is injective.  The condition is best stated in the language of cube complexes, that is, cell complexes whose cells are all cubes.  The link of a vertex $v$ in a cube complex $X$ can be given the structure of a simplicial complex with one $d$-simplex for each $(d+1)$-cube of $X$ incident to $v$.

There are two cube complexes relevant to us.  First, any $\Gamma$-dissection $\{c_i\}$ of a closed surface $\Sigma_h$ gives rise to a cell structure on $\Sigma_h$ where each intersection point of the $c_i$ is a vertex; the dual cell decomposition is a cube complex $X_{\{c_i\}}$ with one square for each intersection point of the $c_i$.

There is also a cube complex $X_\Gamma$, sometimes called the Salvetti complex for $A(\Gamma)$, which is a natural $K(A(\Gamma),1)$ space.  It has one vertex, one loop for each vertex of $\Gamma$, one 2-torus for each edge of $\Gamma$ (glued along the corresponding commutator $[v_i,v_j]$), and, more generally, one $k$-torus for each complete graph on $k$ vertices in $\Gamma$.

The label-reading map induces a cellular map $X_{\{c_i\}} \to X_\Gamma$.  We say that this map satisfies the \emph{link condition} if for each vertex $v$ of $X_{\{c_i\}}$, the induced map from the link of $v$ to the link of the vertex of $X_\Gamma$ is injective and the image of the link of $v$ is a full subcomplex of the link of the vertex of $X_\Gamma$.

The following is a special case of a theorem of Crisp--Wiest \cite{CW}.

\begin{theorem}
\label{theorem:crisp wiest}
Let $h \geq 2$, let $\Gamma$ be a graph, and let $\{c_i\}$ be a $\Gamma$-dissection of $\Sigma_h$.  If the induced map $X_{\{c_i\}} \to X_\Gamma$ satisfies the link condition, then the  label-reading map $\pi_1(\Sigma_h) \to A(\Gamma)$ is injective.
\end{theorem}

\subsection{Right-angled Artin groups into right-angled Artin groups} We next give explicit embeddings of right-angled Artin groups into right-angled Artin groups due to Kim \cite{K}. 

The \emph{opposite} of a graph $\Gamma$ is the graph $\overline \Gamma$ with the same vertex set as $\Gamma$ and with the property that two vertices of $\overline \Gamma$ are connected by an edge if and only if they are not connected by an edge in $\Gamma$.

The basic idea of Kim's embeddings is that, if we collapse a connected subgraph of $\overline \Gamma$ in order to obtain a graph $\overline \Gamma'$, then there is often a natural injective homomorphism $A(\Gamma') \to A(\Gamma)$.

Let $\Gamma$ be a finite graph, and $S$ a subset of its vertex set $\{v_i\}$.  The induced subgraph of $\Gamma$ corresponding to $S$ is the graph whose vertex set is $S$ and whose edge set is the set of edges of $\Gamma$ with both endpoints in $S$.  Say that the subset $S$ is \emph{anti-connected} in $\Gamma$ if the induced subgraph of $S$ in $\overline \Gamma$ is connected.

The \emph{contraction} $CO(\Gamma,S)$ is the graph obtained by collapsing the induced subgraph of $\Gamma$ with respect to $S$ and removing any edge-loops and repeated edges.  The \emph{co-contraction} of $\Gamma$ with respect to $S$ is:
\[ \overline{CO}(\Gamma,S) \cong \overline{CO\left (\overline \Gamma,S\right )} \]
We denote the vertices of $\overline{CO}(\Gamma,S)$ by $\{v_i\, |\, i \notin S\} \cup \{v_S \}$.

If $w$ is any word in the $\{v_i^{\pm 1} \,|\, i \in S\}$, then there is a homomorphism $A(\overline{CO}(\Gamma,S)) \to A(\Gamma)$ given by
\[ v_i \mapsto \begin{cases} v_i & v_i \neq v_S \\ w & v_i = v_S. \end{cases} \]
That this is a homomorphism follows immediately from the defining presentations and the definition of co-contraction.

The only question that remains: how can we choose $w$ so that the homomorphism is injective?  Kim gives the following sufficient condition.  Write 
\[ w = v_{i_1}^{\epsilon_1} \cdots v_{i_m}^{\epsilon_m} \]
where $\epsilon_i \in \{1,-1\}$ for all $i$.  If we delete each $\epsilon_i$ and some subset of the $v_{i_j}$ from the expression for $w$, we obtain a sequence of vertices in $\Gamma$.  We say that the sequence of vertices \emph{appears in $w$}.

Kim proved the following theorem \cite{K}; see \cite{B} and \cite{KK} for other proofs.

\begin{theorem}
\label{theorem:kim}
Let $\Gamma$ be a finite graph with vertex set $\{v_i\}$, and let $S$ be an anti-connected subset of $\{v_i\}$.  Let $w$ be a word in the set $\{v_i^{\pm 1} | v_i \in S \}$ with the property that, for each ordered pair of vertices in $S$, there is an edge path in the induced graph $\overline \Gamma_S$ from one vertex to the other so that the corresponding sequence of vertices appears in $w$.  Then the homomorphism
\[ v_i \mapsto \begin{cases} v_i & v_i \neq v_S \\ w & v_i = v_S \end{cases} \]
is injective.
\end{theorem}

\subsection{Right-angled Artin groups into the braid group}  There are a number of known methods for constructing embeddings of right-angled Artin groups into braid groups (see the end of Section~\ref{sec:main proof}).  Here, we will present a family of embeddings (for a specific type of Artin group) due to L\"onne.  

The Birman--Ko--Lee \cite{BKL} generating set for the braid group $B_n$ has one generator $\beta_{i,j}$ for every pair $1 \leq i < j \leq n$.  These are square roots of the standard generators for the pure braid group.  Also, the generators $\beta_{i,i+1}$ are precisely the usual generators $\sigma_i$ for $B_n$.  The generator $\beta_{1,n}$ is equal to $\left(\sigma_2 \cdots \sigma_{n-1}\right)^{-1}\sigma_1\left(\sigma_2 \cdots \sigma_{n-1}\right)$.

Let $M = (m_{i,j})$ be a symmetric, nonnegative integer matrix, and let $B_n^M$ denote the subgroup of $B_n$ generated by the $\beta_{i,j}^{m_{i,j}}$.  L\"onne gives the following description of this subgroup \cite{L}.

\begin{theorem}
\label{theorem:lonne}
Let $M = (m_{i,j})$ be a symmetric, nonnegative integer matrix where each $m_{i,j}$ is not equal to 1 or 2.  The group $B_n^M$ is a right-angled Artin group with the presentation
\begin{align*}
 B_n^M \cong \langle \beta_{i,j}^{m_{i,j}} \text{ with } m_{i,j} \neq 0\ |\  [\beta_{i,j}^{m_{i,j}} , \beta_{k,l}^{m_{k,l}} ] = 1 \text{ when } & i < j < k < l \\ \text{ or } & i < k < l < j \rangle
\end{align*}
\end{theorem}

\subsection{The braid group to the mapping class group}
\label{sec:braid to mcg}

\begin{figure}
\labellist
\small\hair 2pt
\pinlabel {$c_1$} [] at 18 52 
\pinlabel {$c_2$} [] at 45 52 
\pinlabel {$c_3$} [] at 72 52 
\pinlabel {$c_4$} [] at 99 52 
\pinlabel {$c_5$} [] at 126 52 
\pinlabel {$c_6$} [] at 153 52 
\pinlabel {$c_{2g}$} [] at 228 52 
\endlabellist
\centering \includegraphics[scale=1.25]{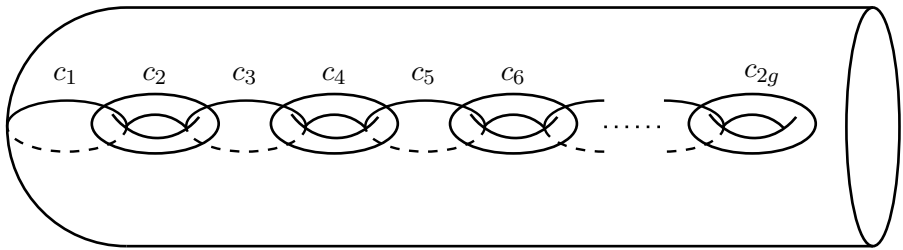}
\caption{The standard embedding of $B_{2g+1}$ in $\MCG(\Sigma_g^1)$.}
\label{figure:braid}
\end{figure}

There are two natural homomorphisms\footnote{Because in the braid group we traditionally compose elements left to right and in the mapping class group we compose right to left, the first map here is really an anti-homomorphism.  This is the only place where we switch the order of composition, and so our maps $\pi_1(\Sigma_h) \to \MCG(\Sigma_g)$ will also be anti-homomorphisms.  We could remedy the situation by composing elements of $\pi_1(\Sigma_h)$, $A(\overline C_{2g+1})$, and $B_{2g+1}$ from right to left.  Instead, we abuse notation by referring to our anti-homomorphisms as homomorphisms.}
\[ B_{2g+1} \to \MCG(\Sigma_g^1) \to \MCG(\Sigma_g). \]
The first map is defined by the rule $\sigma_i \mapsto T_{c_i}$, where $\{\sigma_i\}$ is the standard generating set for $B_{2g+1}$  and the $c_i$ are the curves in $\Sigma_g^1$ shown in Figure~\ref{figure:braid}.  This map is a homomorphism because Dehn twists about curves intersecting once satisfy the braid relation \cite[Proposition 3.11]{FM}.  We have the following theorem of Birman--Hilden; see \cite[Theorem 9.2]{FM}.

\begin{theorem}
\label{thm:BH}
The homomorphism $B_{2g+1} \to \MCG(\Sigma_g^1)$ given by $\sigma_i \mapsto T_{c_i}$ is injective.
\end{theorem}

The image of $B_{2g+1}$ in $\MCG(\Sigma_g^1)$ is precisely the hyperelliptic mapping class group of $\MCG(\Sigma_g^1)$.  This is the subgroup of $\MCG(\Sigma_g^1)$ consisting of elements that have a representative commuting with a particular hyperelliptic involution of $\Sigma_g^1$; see \cite[Section 9.4]{FM}.  Because our monodromy homomorphisms factor through this map, all of our bundles will be hyperelliptic surface bundles over surfaces.

The second map, $\MCG(\Sigma_g^1) \to \MCG(\Sigma_g)$,  is obtained by gluing a disk to the boundary of $\Sigma_g^1$.  Each homeomorphism of $\Sigma_g^1$ can be extended by the identity to a homeomorphism of $\Sigma_g$, and this induces a well-defined map on the level of mapping class groups \cite[Theorem 3.18]{FM}.

Under this homomorphism,  the Dehn twist about $c_i$ maps to the Dehn twist about the image of this curve under the inclusion $\Sigma_g^1 \to \Sigma_g$.  The kernel is isomorphic to the fundamental group of the unit tangent bundle of $\Sigma_g$ \cite[Section 4.2.5]{FM}.  However, the intersection of the kernel of this map with the image of $B_{2g+1}$ is equal to $Z(B_{2g+1}) \cong \Z$; see \cite[Theorem 3.1]{BM}.

\section{Constructing indecomposable surface bundles} \label{Sec: Construction}

In this section we prove the Main Theorem, that is, for any $g \geq 2$ and $h \geq 1$ we explicitly construct infinitely many genus $g$ surface bundles over genus $h$ surfaces that are both fiber sum indecomposable and section sum indecomposable.  The examples are all distinct in that they are pairwise homotopy inequivalent.

To begin, in Section~\ref{sec:construction} we explain one specific construction of indecomposable bundles over surfaces of genus at least 2.  Then, in Section~\ref{sec:main proof} we complete the proof the Main Theorem using the construction from Section~\ref{sec:construction} (the case there the base surface is a torus is dealt with separately). At the end, we discuss various geometric properties of our bundles, namely, holomorphicity, symplecticity, and signature.

\subsection{A construction of indecomposable bundles}
\label{sec:construction}


We will now give an explicit construction of surface bundles over surfaces that are both fiber sum indecomposable and section sum indecomposable.  We will give one bundle $X_n=X_n(g,h)$ for each $g \geq 2$, $h \geq 1$, and $n \geq 3$.  Recall that, by the Fiber Sum Criterion and the Section Sum Criterion, it suffices to construct injective, irreducible homomorphisms:
\[ \pi_1(\Sigma_h) \to \MCG(\Sigma_g) \]
for varying $g$ and $h$ (and $n$).  

Let $C_k$ denote the $k$-cycle graph.  Our recipe for explicit irreducible embeddings $\pi_1(\Sigma_h) \to \MCG(\Sigma_g)$ (for any $g,h \geq 2$) is broken into four steps:
\[ \pi_1(\Sigma_h) \to A(\overline C_5)  \to A(\overline C_{2g+1}) \to B_{2g+1} \to \MCG(\Sigma_g). \]
At each stage, we give explicit maps.  After giving the construction, we will check that the composition $\pi_1(\Sigma_h) \to \MCG(\Sigma_g)$ is injective and irreducible.  


\fbox{Step 1}\ \ $\pi_1(\Sigma_h) \to A(\overline C_5)$.

To fix notation set
\[ A(\overline C_5) = \langle v_1,\dots,v_5\, | \, [v_1,v_3]=[v_2,v_4]=[v_3,v_5]=[v_4,v_1]=[v_5,v_2]=1  \rangle. \]
Note that $\overline C_5 \cong C_5$, and so $A(\overline C_5) \cong A(C_5)$.  We use the former because it will simplify the notation later.

\begin{figure}
\labellist
\small\hair 2pt
\pinlabel {$v_2$} [] at -11 80 
\pinlabel {$v_1$} [] at 105 33 
\pinlabel {$v_4$} [] at 53 126 
\pinlabel {$v_5$} [] at 53 103 
\pinlabel {$v_5$} [] at 200 103 
\pinlabel {$v_4$} [] at 75 19 
\pinlabel {$v_1$} [] at 125 110 
\pinlabel {$v_3$} [] at 155 102 
\pinlabel {$v_2$} [] at 228 80 
\endlabellist
\centering \includegraphics{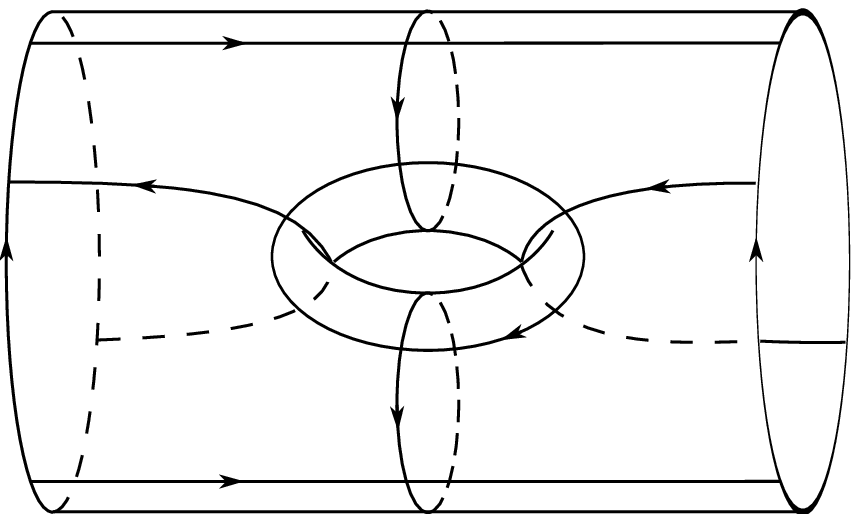}
\caption{Each curve labeled $v_i$ corresponds to the generator $v_i$ of $A(\overline C_5)$.}
\label{figure:crispwiest}
\end{figure}

Consider the surface with boundary shown in Figure~\ref{figure:crispwiest}.  We have drawn several oriented simple closed curves and arcs there.  We obtain a collection of oriented simple closed curves in $\Sigma_2$ by gluing together the two boundary components.  This collection of curves is a $\overline C_5$-dissection of $\Sigma_2$ which appears in the paper of Crisp--Wiest \cite[Figure 2]{CW}.

More generally, we obtain a $\overline C_5$-dissection of $\Sigma_h$ for any $h \geq 2$ by stacking $h-1$ copies of the surface with boundary in Figure~\ref{figure:crispwiest} end to end and then gluing the two resulting boundary components.  The labels are inherited directly from the original curves.  (Equivalently, the dissection of $\Sigma_h$ is the one induced by the $(h-1)$-fold cover of $\Sigma_2$ dual to $v_2$.)

As in Section~\ref{sec:crisp wiest}, the $\overline C_5$-dissection of $\Sigma_h$ gives rise to a label-reading homomorphism $\pi_1(\Sigma_h) \to A(\overline C_5)$.  We would like for this homomorphism to be injective.  It is straightforward to check the link condition of Theorem~\ref{theorem:crisp wiest} (see \cite[page 451]{CW} for details), and so the label-reading map $\pi_1(\Sigma_h) \to A(\overline C_5)$ is indeed injective, as desired.

\begin{figure}
\labellist
\small\hair 2pt
\pinlabel {$\gamma_1$} [] at 67 47 
\pinlabel {$\gamma_2$} [] at 195 103 
\pinlabel {$\delta_1$} [] at 177 115 
\pinlabel {$\delta_2$} [] at 67 80 
\endlabellist
\centering \includegraphics{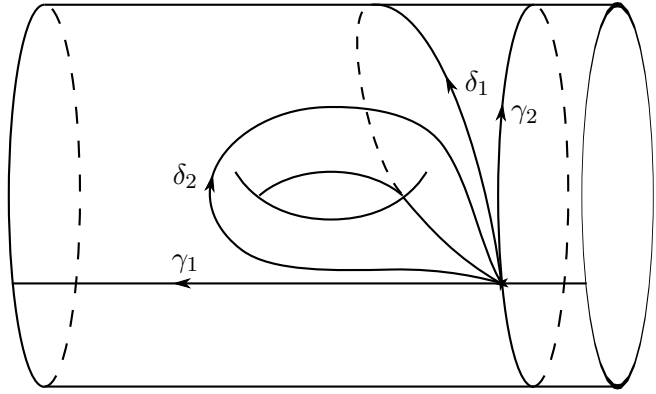}
\caption{A standard set of generators for $\pi_1(\Sigma_2)$.}
\label{figure:generators}
\end{figure}

If we utilize the generating set for $\pi_1(\Sigma_2)$ shown in Figure~\ref{figure:generators}, the label-reading map $\pi_1(\Sigma_2) \to A(\overline C_5)$ has the following effect:
\begin{align*}
\gamma_1 &\mapsto v_1^{-1}v_2 \\
\gamma_2 &\mapsto v_5^{-1}v_4v_5v_4 \\
\delta_1 &\mapsto v_5^{-1} v_4 v_5 v_3 \\
\delta_2 &\mapsto v_1^{-1}v_5^{-1} v_1 v_5
\end{align*}

\begin{figure}
\labellist
\small\hair 2pt
\pinlabel {$\delta_4$} [] at 39 47 
\pinlabel {$\gamma_1$} [] at 45 16 
\pinlabel {$\delta_1$} [] at 223 80 
\pinlabel {$\gamma_2$} [] at 245 80 
\pinlabel {$\delta_2$} [] at 174 48 
\pinlabel {$\delta_3$} [] at 90 80 
\endlabellist
\centering \includegraphics{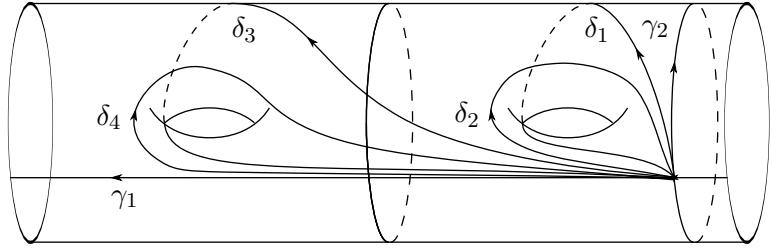}
\caption{A set of generators for $\pi_1(\Sigma_3)$ adapted to the covering $\Sigma_3 \to \Sigma_2$ dual to $v_2$.}
\label{figure:generators2}
\end{figure}

Since our dissection of $\Sigma_h$ is induced by the $(h-1)$-fold cover $\Sigma_h \to \Sigma_2$ dual to $v_2$, we can write our map $\pi_1(\Sigma_h) \to A(\overline C_5)$ as the composition $\pi_1(\Sigma_h) \to \pi_1(\Sigma_2) \to A(\overline C_5)$.  Consider the generating set for $\pi_1(\Sigma_h)$ suggested by Figure~\ref{figure:generators2}.  With respect to these generators, the map $\pi_1(\Sigma_h) \to \pi_1(\Sigma_2)$ is given by:
\begin{align*}
\gamma_1  & \mapsto \gamma_1^{h-1} \\
\gamma_2 & \mapsto \gamma_2
\end{align*}
and
\begin{align*}
\delta_{2k+1} & \mapsto \gamma_1^{k} \delta_1 \gamma_1^{-k} \\
\delta_{2k+2} & \mapsto \gamma_1^{k} \delta_2 \gamma_1^{-k} \\
\end{align*}
for $0 \leq k \leq h-2$.

\fbox{Step 2}\ \ $A(\overline C_5) \to A(\overline C_{2g+1})$.

Denote the vertices of $C_{2g+1}$ by $v_1,\dots,v_{2g+1}$ (in cyclic order), and let $S=\{v_4, \dots, v_{2g}\}$.  We have
\[ \overline{CO}(\overline C_{2g+1},S) = \overline{CO\left(C_{2g+1},S\right)} \cong \overline C_5.\]
The group $A(\overline{CO}(\overline C_{2g+1},S))$ is generated by $v_1$, $v_2$, $v_3$, $v_S$, and $v_{2g+1}$.  There is an isomorphism $A(\overline C_5) \to A(\overline{CO}(\overline C_{2g+1},S))$ defined by \[ (v_1,v_2,v_3,v_4,v_5) \mapsto (v_1,v_2,v_3,v_S,v_{2g+1}).\]

Let 
\[ w = v_4 v_5^{-1} \cdots v_{2g-1}^{-1} v_{2g} v_{2g-1}^{-1} \cdots v_5^{-1} v_4 \in A(\overline C_{2g+1}). \]
Notice that the sequences $v_4,\dots,v_{2g}$ and $v_{2g},\dots,v_4$ both appear in $w$.  By Theorem~\ref{theorem:kim}, the homomorphism 
\[ A(\overline C_5) \stackrel{\cong}{\to} A(\overline{CO}(\overline C_{2g+1},S)) \to A(\overline C_{2g+1}) \]
defined by
\[ v_i \mapsto \begin{cases} v_i & i=1,2,3 \\ w & i =4 \\ v_{2g+1} & i=5. \end{cases} \]
is injective.


\fbox{Step 3}\ \ $A(\overline C_{2g+1}) \to B_{2g+1}$

Let $M = (m_{i,j})$ be given by: 
\[ m_{i,j} = \begin{cases} n & j \equiv i \pm 1 {\mod 2g+1} \\ 0 & \text{otherwise}   \end{cases} \]
By Theorem~\ref{theorem:lonne}, the group $B_{2g+1}^M$ has the presentation
\begin{align*} \langle \sigma_1,\dots,\sigma_{2g},\beta_{1,2g+1}\, |\, & [\sigma_i,\sigma_j]=1 \text{ for } |i-j| > 1,\\& \ \ \ [\beta_{1,2g+1},\sigma_i]=1 \text{ for } i \notin \{1,2g\} \rangle. \end{align*}
There is thus an isomorphism
\[ A(\overline C_{2g+1}) \stackrel{\cong}{\to} B_{2g+1}^M < B_{2g+1} \]
given by
\[ v_i \mapsto \begin{cases} \sigma_i & 1 \leq i \leq 2g \\ \beta_{1,2g+1} & i = 2g+1. \end{cases} \]



\fbox{Step 4}\ \ $B_{2g+1} \to \MCG(\Sigma_g)$

\begin{figure}
\labellist
\small\hair 2pt
\pinlabel {$c_{2g+1}$} [] at 67 8 
\endlabellist
\centering \includegraphics[scale=1.25]{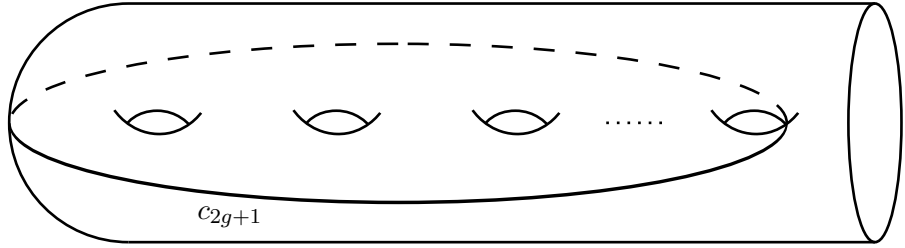}
\caption{The simple closed curve $c_{2g+1}$.}
\label{figure:braid2}
\end{figure}

The map here was already described in Section~\ref{sec:braid to mcg}.  We will make several remarks.  First, the image of $\beta_{1,2g+1} \in B_{2g+1}$ is the Dehn twist $T_{c_{2g+1}}$, where $c_{2g+1}$ is the simple closed curve shown in Figure~\ref{figure:braid2}.   In terms of the $T_{c_i}$, we can write
\[ T_{c_{2g+1}} = \left(T_{c_{2g}}\cdots T_{c_3}T_{c_2}\right)T_{c_1}\left(T_{c_{2g}}\cdots T_{c_3}T_{c_2}\right)^{-1}. \]
Also, the class of $c_{2g+1}$ in $H_1(\Sigma_g^1)$ is $[c_1]-[c_2]+[c_3] - \cdots + [c_{2g-1}] - [c_{2g}]$.

Finally, for large $n$, the injectivity of our map $A(\overline C_{2g+1}) \to \MCG(\Sigma_g)$ can also be deduced directly from a theorem of Koberda \cite[Theorem 1.1]{Ko}.  We have taken the slightly more circuitous route in order to gain the explicitness.  


\fbox{Step 5}\ \ \emph{Injectivity}

Our homomorphism $\pi_1(\Sigma_h) \to \MCG(\Sigma_g)$ is the composition of the above homomorphisms:
\[ \pi_1(\Sigma_h) \to A(\overline C_{2g+1},S)  \to A(\overline C_{2g+1}) \to B_{2g+1} \to \MCG(\Sigma_g^1) \to \MCG(\Sigma_g). \]
The only map in this sequence that is not injective is the last one.  Recall that the kernel of the composition
\[ B_{2g+1} \to \MCG(\Sigma_g^1) \to \MCG(\Sigma_g) \]
is $Z(B_{2g+1})$.  Under an injective homomorphism, noncentral elements cannot map to central elements.  As $\pi_1(\Sigma_g)$ has trivial center for \linebreak $g \geq 2$, it follows that the composition $\pi_1(\Sigma_h) \to \MCG(\Sigma_g)$ is injective, as desired.


\fbox{Step 6}\ \ \emph{Irreducibility}

The Nielsen--Thurston classification theorem states that each element of $\MCG(\Sigma_g)$ falls into one of three categories: periodic, reducible, or pseudo-Anosov; see \cite[Theorem 13.2]{FM}.  Moreover, pseudo-Anosov elements are neither periodic nor reducible.

In order to prove that our map $\pi_1(\Sigma_h) \to \MCG(\Sigma_g)$ is irreducible, we will use Penner's construction of pseudo-Anosov elements of $\MCG(\Sigma_g)$, as follows \cite{P}.  In the statement, we say that two simple closed curves in a surface are in minimal position if they intersect the minimal number of times with respect to their two homotopy classes.

\begin{theorem}
\label{theorem:penner}
Let $A=\{a_i\}$ and $B=\{b_i\}$ be two collections of essential simple closed curves in $\Sigma_g$ with the following properties:
\begin{enumerate}
\item We have $a_i \cap a_j = \emptyset$ and $b_i \cap b_j = \emptyset$ for any choices of $i$ and $j$.
\item For any $i$ and $j$, the curves $a_i$ and $b_j$ are in minimal position.
\item The complement of $A \cup B$ in $\Sigma$ is a union of disks.
\end{enumerate}
Then any product of the $T_{a_i}$ and $T_{b_j}^{-1}$ where each $a_i$ and each $b_j$ appears at least once represents a pseudo-Anosov element of $\MCG(\Sigma_g)$.
\end{theorem}

In our situation, the elements $\gamma_1$, $\gamma_2$, and $\delta_1$ in $\pi_1(\Sigma_h)$ shown in Figure~\ref{figure:generators2} map to $(v_1^{-1}v_2)^{h-1}$,  $v_5^{-1}v_4v_5v_4$, and $v_5^{-1}v_4v_5v_3$ in $A(\overline C_5)$.  Thus the product $\gamma_2^{-1}\delta_1\gamma_1^{-1}$ maps to
\[ (v_5^{-1}v_4v_5v_4)^{-1}(v_5^{-1}v_4v_5v_3)(v_1^{-1}v_2)^{-1} = v_4^{-1}v_3(v_2^{-1}v_1)^{h-1} \]
in $A(\overline C_5)$, and then
\[ w^{-1}v_3(v_2^{-1}v_1)^{h-1} = \left(v_4^{-1} v_5 \cdots v_{2g-1} v_{2g}^{-1} v_{2g-1} \cdots v_5 v_4^{-1} \right) v_3(v_2^{-1}v_1)^{h-1} \]
in $A(\overline C_{2g+1})$. 
The image in $\MCG(\Sigma_g^1)$ is:
\[ \left(T_{c_1}^{n}T_{c_2}^{-n}\right)^{h-1}\left(T_{c_3}^n\right)\left(T_{c_4}^{-n} T_{c_5}^{n} \cdots T_{c_{2g-1}}^{n} T_{c_{2g}}^{-n} T_{c_{2g-1}}^{n} \cdots T_{c_5}^{n}T_{c_4}^{-n}\right), \]
where the $c_i$ are as shown in Figure~\ref{figure:braid}.  It follows from Theorem~\ref{theorem:penner} that this maps to a pseudo-Anosov element of $\MCG(\Sigma_g)$ (set $A=\{c_i\, |\, i \text{ odd} \}$ and $B=\{c_i \,|\, i \text{ even}\}$).

The irreducibility of our monodromy follows from general principles, namely, the classification of subgroups of the mapping class group; see \cite{BLM,ivanov,Ma,Mo}.  For our purposes, it is simpler to simply produce an explicit pseudo-Anosov monodromy using Penner's theorem than to explain the more general theory.


\fbox{Step 7}\ \ \emph{Monodromy factorization}

We first deal with the case $h=2$, since it is notationally simpler.  In Figure~\ref{figure:generators} we have drawn generators for $\pi_1(\Sigma_2)$.  Referring to Figure~\ref{figure:crispwiest}, we see that the label-reading map $\pi_1(\Sigma_2) \to A(\overline C_5)$ has the following effect:
\[ (\gamma_1,\gamma_2,\delta_1,\delta_2) \mapsto (v_1^{-1}v_2\ ,\   v_5^{-1}v_4v_5v_4\ ,\   v_5^{-1} v_4 v_5 v_3\ ,\  v_1^{-1}v_5^{-1} v_1 v_5). \]
It follows that the images of the generators for $\pi_1(\Sigma_2)$ in $\MCG(\Sigma_g)$ are:
\begin{align*}
\gamma_1  & \mapsto T_{c_2}^{n}T_{c_1}^{-n} \\
\gamma_2 & \mapsto W T_{c_5}^{n} W T_{c_5}^{-n} \\
\delta_1 & \mapsto T_{c_3}^nT_{c_5}^{n}WT_{c_5}^{-1} \\
\delta_2 & \mapsto T_{c_5}^nT_{c_1}^{n}T_{c_5}^{-n}T_{c_1}^{-n}
\end{align*}
where $W = T_{c_4}^{n} T_{c_5}^{-n} \cdots T_{c_{2g-1}}^{-n} T_{c_{2g}}^{n} T_{c_{2g-1}}^{-n} \cdots T_{c_5}^{-n}T_{c_4}^{n}$.

We conclude that there is a fiber sum and section sum indecomposable genus $g \geq 2$ surface bundle over a genus $h=2$ surface prescribed by the monodromy factorization
\[
 [
 T_{c_2}^{-n}T_{c_1}^{n}
  \, , 
  W T_{c_5}^{n} W T_{c_5}^{-n}
 ]
\,
[
T_{c_3}^nT_{c_5}^{n}WT_{c_5}^{-1}
\, , 
T_{c_5}^nT_{c_1}^{-n}T_{c_5}^{-n}T_{c_1}^{n}
] 
 =1.
\]
As above, our bundles over genus $h \geq 3$ surfaces are obtained by pulling back the above bundles via the $(h-1)$-fold cyclic cover $\Sigma_h \to \Sigma_2$ dual to $v_2$.  In Step 1, we wrote down the map $\pi_1(\Sigma_h) \to \pi_1(\Sigma_2)$ in terms of standard generators for $\pi_1(\Sigma_h)$.  From here we can easily write down the monodromy factorization. 
\subsection{Infinite families}
\label{sec:main proof}

In this section we complete the proof of the Main Theorem.  The case $h=1$ will be handled from scratch.  The case $h \geq 2$ will be dealt with by showing that the $X_n(g,h)$ constructed in Section~\ref{sec:construction} are all distinct (for $n$ prime).

\bigskip {\bf Computing first homology.} Before delving into the proof, we recall how the first homology of a bundle is computed.  Let $X \to B$ be a fiber bundle with fiber $F$.  Via the monodromy, the group $\pi_1(B)$ acts on $H_1(F)$. Using this action we obtain
\[ H_1(X) \cong H_1(B) \oplus \left( H_1(F)/\pi_1(B) \right). \]
We can compute this quotient by fixing generators for $\pi_1(B)$ and $H_1(F)$.  There is then one (possibly trivial) relation for the action of each generator of $\pi_1(B)$ on each generator of $H_1(F)$.

As all of our monodromies are expressed in terms of Dehn twists, it is enough to repeatedly use the formula:
\[ T_b^k([a]) = [a]-k\, \widehat i (a,b)[b] \]
for oriented simple closed curves $a$ and $b$; here $\widehat i(a,b)$ denotes algebraic intersection.  We will do one sample calculation of this sort in the proof below and leave the rest of the (straightforward) calculations to the reader.

\bigskip {\bf Bundles over tori.} We first deal with the case $h=1$.  Since the torus cannot be written as a nontrivial connect sum, fiber sum decomposability is immediately ruled out.  By the Section Sum Criterion, it suffices to produce bundles whose monodromies $\pi_1(T^2) \to \MCG(\Sigma_g)$ are irreducible.

Let $\phi$ be any irreducible element of $\MCG(\Sigma_g)$.  The monodromy factorization
\[ [\phi, 1] =1 \]
corresponds to a section sum indecomposable bundle, as desired. (Observe that this construction amounts to taking the product of $S^1$ and a $\Sigma_g$ bundle over $S^1$ with irreducible monodromy.)

In order to obtain infinitely many bundles whose total spaces are distinct up to homotopy equivalence, it suffices to find infinitely many such $\phi$ whose actions on $H_1(F)$ are distinct in the sense that the groups $H_1(F)/\pi_1(B)$ are nonisomorphic.  One example of such an infinite family is:
\[ \phi_k = T_{c_1}^k T_{c_2}^{-1} T_{c_3} T_{c_4}^{-1} \cdots T_{c_{2g-1}} T_{c_{2g}}^{-1}. \]
Each $\phi_k$ with $k > 0$ is pseudo-Anosov (hence irreducible) by Penner's Theorem~\ref{theorem:penner}. (Again set $A=\{c_i\, |\, i \text{ odd} \}$ and $B=\{c_i \,|\, i \text{ even}\}$). Denote by $X_k$ the bundle over $T^2$ with monodromy factorization $[\phi_k, 1] =1$. 

Choose orientations for the $c_i$ so that the algebraic intersection numbers $\widehat i(c_i,c_{i+1})$ are equal to 1, and denote by $[c_i]$ the resulting elements of $H_1(F) \subseteq H_1(X_k)$.  As discussed above, we can compute $H_1(X_k)$ by adding relations $\phi_k([c_i])=[c_i]$ for each $i$.  In the case $i=1$, we obtain: 
\begin{align*}
[c_1] &= \phi_k([c_1]) \\
 &= T_{c_1}^k T_{c_2}^{-1} T_{c_3} T_{c_4}^{-1} \cdots T_{c_{2g-1}} T_{c_{2g}}^{-1}([c_1]) \\
 &= T_{c_1}^k T_{c_2}^{-1} ([c_1]) \\
 &= T_{c_1}^k([c_1] +[c_2]) \\
 &= [c_1]+([c_2]+k[c_1]) \\
\end{align*}
So $k[c_1]=-[c_2]$.  There are $2g-1$ similar calculations.  They yield the relations $-[c_2]=[c_3], [c_3]=-[c_4], \dots, [c_{2g-1}]=-[c_{2g}]$, and $[c_{2g}]=0$.  We thus conclude that 
\[ H_1(X_k) \cong  H_1(T^2) \oplus \Z/k\Z \cong \Z^2 \oplus \Z/k\Z.\]
In particular, the total spaces of the $X_k$ are pairwise homotopy inequivalent, as desired.

\bigskip{\bf Bundles over higher genus surfaces.} Now let $h \geq 2$.  For any $g \geq 2$ and $n \geq 3$, let $X_n=X_n(g,h)$ be the bundle constructed in Section~\ref{sec:construction}. Recall that all these bundles are fiber sum and section sum indecomposable. 
We will show that the family $\{X_n : n \text{ prime}\}$ consists of $4$-manifolds with distinct first homology groups, and thus of $4$-manifolds that are pairwise homotopy inequivalent. 

Fix some $X_n$ and denote its base by $B$ and its fiber by $F$. We have
\[ H_1(X_n) \cong \Z^{2h} \oplus \left( H_1(F)/\pi_1(B) \right). \]
where $\pi_1(B)$ acts via the monodromy. As above, we take the $[c_i]$ as a basis for $H_1(F)$. 

For our purposes it will suffice to make two observations.

Firstly, since the monodromy action of $\pi_1(B)$ on $H_1(F; \Z/n\Z)$ is trivial, we easily calculate
\[ H_1(X_n ; \Z/n\Z) \cong (\Z/n\Z)^{2h} \bigoplus (\Z/n\Z)^{2g} = \bigoplus (\Z/n\Z)^{2h+2g} .\]
It then follows from the universal coefficients theorem that, for $n=p$ prime, $H_1(X_n)=H_1(X_n; \Z)$ has a total of $2h+2g$ direct summands, each isomorphic to either $\Z$ or $\Z/p^k\Z$ for some $k \geq 1$. 

Secondly, when we calculate $H_1(X_n)$ as above, we find many relations between the $2h+2g$ generators $c_1, \ldots, c_{2g}$. For example, it is straightforward to check that the action of $\gamma_1$ on $[c_1]$ gives a nontrivial relation (in the $h=2$ case the relation is $n \, [c_2] = 0$).

Combining the two observations, we see that at least one summand of $H_1(X_p;\Z)$ must be of the form $\Z/p^k\Z$.  It follows that the $H_1(X_p;\Z)$ are all distinct, and the $\{X_p : p \text{ prime} \}$ are pairwise homotopy inequivalent, as desired.  This completes the proof of our Main Theorem.


\bigskip{\bf Variations.} While our construction of surface bundles is general enough to give the infiniteness result of the Main Theorem, we would like to point out some of the ways in which the construction can be altered in order to give other examples of indecomposable surface bundles over surfaces.

First, there are many $C_5$-dissections of $\Sigma_h$.  There are also completely different approaches to embedding surface groups into right-angled Artin groups; see the work of Servatius--Droms--Servatius \cite{SDS}, Crisp--Wiest \cite{CW}, R\"over \cite{Ro}, Crisp--Sageev--Sapir \cite{CSS}, Kim \cite{K}, and Bell \cite{B}.

Next, there are other known methods for embedding right-angled Artin groups into other right-angled Artin groups; see the papers of Bestvina--Kleiner--Sageev \cite{BKS}, Hsu--Wise \cite{HW}, and Green \cite{G}.

Various embeddings of right-angled Artin groups into braid groups and mapping class groups are given by Collins \cite{Co}, Humphries \cite{Hu}, Crisp--Paris \cite{CP}, Crisp--Wiest \cite{CW}, Clay--Leininger--Mangahas \cite{CLM}, and Koberda \cite{Ko}.

Finally, there are other approaches to embedding surface groups into mapping class groups that do not pass through the theory of right-angled Artin groups; see, for instance, the works of Gonz{\'a}lez-D{\'{\i}}ez--Harvey \cite{GDH}, and Leininger--Reid \cite{LR}.

Because of all of this work, there are many possible paths for embedding surface groups into mapping class groups.  We chose to focus on the constructions of Crisp--Wiest, Kim, L\"onne, and Birman--Hilden since this approach is both completely explicit and strikingly simple.


\bigskip{\bf Final remarks.} A surface bundle over a surface $(X,f)$ is called \emph{symplectic} if $X$ admits a symplectic form for which all the fibers are symplectic subsurfaces. It is well-known that when the surface bundles $(X_i, f_i)$, and the sections $S_i$ when involved, are symplectic, then the fiber sum and section sum operations can be performed symplectically \cite{GS}. One can thus ask: when does a symplectic surface bundle over a surface fail to decompose into symplectic surface bundles of smaller fiber or base genera?

It follows from a classical argument of Thurston \cite{Thurston} that if the fiber $F$ of a surface bundle $(X,f)$ is nonzero in $H_2(X ; \R)$, then there is a symplectic form on $X$ so that $F$ is symplectic.  The symplectic form can be chosen so that both the fiber and a prescribed section are symplectic subsurfaces in $X$. The first chern class of an almost complex structure associated to the fibration $f$ gives a class in  $H^2(X; \R)$, which evaluates on $[F]$ as $\chi(F)=2-2g$, which implies that the above homological condition is satisfied whenever $g \neq 1$. We therefore conclude that the surface bundles we construct in this article are all symplectic bundles that cannot be decomposed as fiber sums or section sums of (symplectic) bundles with smaller fiber or base genera.

It is quite often the case that the surface bundle over a surface obtained by performing fiber sum or section sum is not holomorphic for any choice of complex structures on the total space (with either orientation) or the base; see for instance \cite{Baykur}. One can thus inquire whether our surface bundles owe their indecomposability to being holomorphic. This is however ruled out by Parshin's finiteness result, which states that there are only finitely many holomorphic fibrations with fixed fiber genus $g \geq 2$ and base genus $h \geq 0$ \cite{Parshin}. The very fact that we get infinite families of such bundles therefore imply that an infinite subfamily of them are not holomorphic. 

As per the discussion after Theorem~\ref{thm:BH}, all of our bundles are hyperelliptic surface bundles. Such bundles necessarily have signature zero \cite{CLM76, Ham}. In particular, it follows that our bundles are distinct from the earlier examples discussed in the introduction. More to the point, this means that we cannot apply the same signature obstruction to fiber sum decomposability.  


\bibliographystyle{plain}
\bibliography{IndecomposableSurfaceBundles}

\end{document}